\documentclass[reqno]{amsart}

\usepackage{amssymb,amsxtra,latexsym}

\usepackage[sc]{mathpazo}
\usepackage{mathrsfs}
\usepackage[T1]{fontenc}
\usepackage{textcomp}

\usepackage[cmtip,matrix,arrow,curve]{xy}
\CompileMatrices
\newdir{((}{{\lhook\kern-.1em}}
\newdir{))}{{\rhook\kern-.1em}}
\newdir{ >}{{}*!/-5pt/@{>}}

\usepackage{multicol}

\usepackage{booktabs}

\numberwithin{equation}{section}

\renewcommand\dots{\relax\ifmmode\ldots\else$\,\ldots\,$\fi}

\newcommand\note[1]%
{$^\dagger$\marginpar{\footnotesize{$^\dagger${#1}}}}

\divide\count253 by 60 
\divide\count255 by 60
\multiply\count255 by 60 
\advance\count254 by -\count255 

\def\today{\number\year-\ifnum\month<10
0\fi\number\month-\ifnum\day<10 0\fi\number\day}

\def\hour{\ifnum\count253<10
0\number\count253\else\number\count253\fi}

\def\minute{\ifnum\count254<10
0\number\count254\else\number\count254\fi}

\swapnumbers

\newtheorem*{theorem*}{Theorem}
\newtheorem*{lemma*}{Lemma}

\newtheorem*{proposition*}{Proposition}

\newtheorem*{corollary*}{Corollary}

\theoremstyle{definition}

\newtheorem*{example*}{Example}

\newcommand\lie{\mathfrak}

\newcommand{\g}{\lie{g}}
\newcommand{\h}{\lie{h}}

\renewcommand{\t}{\lie{t}}

\newcommand\bb[1]{{\text{\bf#1}}}

\newcommand\Z{\bb{Z}} 

\newcommand\R{\bb{R}} 
\newcommand\C{\bb{C}}

\newcommand\ca{\mathscr}

\newcommand\func[1]{\operatorname{\mathrm{#1}}}

\newcommand\Aut{\func{Aut}}
\newcommand\Ad{\func{Ad}}
\newcommand\ad{\func{ad}}

\renewcommand\det{\func{det}}

\newcommand\Hom{\func{Hom}}
\newcommand\id{\func{id}}

\renewcommand\index{\func{index}}

\newcommand\Pf{\func{Pf}}

\newcommand\vol{\func{vol}}

\newcommand\group[1]{{\text{\bf#1}}}

\newcommand\Spin{\group{Spin}}

\newcommand\GL{\group{GL}}

\newcommand\G{\group{G}}
\newcommand\T{\group{T}}

\newcommand\bigabs[1]{\bigl\lvert#1\bigr\rvert}
\newcommand\norm[1]{\lVert#1\rVert}

\newcommand\quot[1][\kern.3ex]{/\kern-.7ex/_{\kern-.4ex#1}}
\newcommand\bigquot[1][\,\,]{\big/\kern-.85ex\big/_{\!\!#1}}

\newcommand\powl{[\kern-.3ex[}
\newcommand\powr{]\kern-.3ex]}
\newcommand\bigpowl{\bigl[\kern-.6ex\bigl[}
\newcommand\bigpowr{\bigr]\kern-.6ex\bigr]}

\newcommand\sur{\mathrel{\to\kern-1.8ex\to}}
\newcommand\iso{\mathrel{\hookrightarrow\kern-1.8ex\to}}

\newcommand\longto{\longrightarrow}

\newcommand\longsur{\mathrel{\longrightarrow\kern-1.8ex\to}}

\newcommand\dirac{\textit{\dh}}

\newcommand\e{\group{e}}

\newcommand\Pic{\group{Pic}}
\newcommand\TPic{\group{TPic}}

\newcommand\zerodots%
{\smash{\makebox[0pt]{$_{\lower8pt\hbox{\Large$0$}}%
\ddots^{\lower3pt\hbox{\Large$0$}}$}}}
\newcommand\bigzerodots%
{\smash{\makebox[0pt]{$_{\lower6pt\hbox{\Large$0$}}%
\ddots^{\hbox{\Large$0$}}$}}}

\newcommand\eps{\varepsilon}
\renewcommand\phi{\varphi}

\begin{document}


\title{Hans Duistermaat's contributions to Poisson geometry}

\author{Reyer Sjamaar}

\email{sjamaar@math.cornell.edu}

\address{Department of Mathematics, Cornell University, Ithaca, NY
14853-4201, USA}

\date{2011-10-25}


\begin{abstract}
Hans Duistermaat was scheduled to lecture in the 2010 School on
Poisson Geometry at IMPA, but passed away suddenly.  This is a record
of a talk I gave at the 2010 Conference on Poisson Geometry (the week
after the School) to share some of my memories of him and to give a
brief assessment of his impact on the subject.
\end{abstract}


\maketitle



Johannes Jisse (Hans) Duistermaat (1942--2010) earned his doctorate in
1968 at the University of Utrecht under the direction of Hans
Freudenthal.  After holding a postdoctoral position at the University
of Lund, a professorship at the University of Nijmegen, and a visiting
position at the Courant Institute, he returned to Utrecht in 1975 to
take over Freudenthal's chair after the latter's retirement.  He held
this chair until his own retirement in 2007.  He continued to play a
role in the life of the Utrecht Mathematical Institute and kept up his
mathematical activities until he was struck down in March 2010 by a
case of pneumonia contracted while on chemotherapy for cancer.

I got to know Hans Duistermaat as an undergraduate when I took his
freshman analysis course at Utrecht.  At that time he taught the
course from lecture notes that were in style and content close to
Dieudonn\'e's book \cite{dieudonne;foundations}.  Although this was a
smashing success with some students, I think Hans realized he had to
reduce the potency so as not to leave behind quite so many of us, and
over the years the lecture notes grew into his and Joop Kolk's still
very substantial undergraduate textbook
\cite{duistermaat-kolk;multidimensional-real}.  Anyway, I quickly
became hooked and realized that I wanted one day to become his
graduate student.  This came to pass and I finished in 1990 my thesis
work on a combination of two of Hans' favourite topics, Lie groups and
symplectic geometry.

The online Mathematics Genealogy Project (accessed 19 December 2010)
lists twenty-four PhD students under Duistermaat's name.  What the
website does not tell you is that his true adviser was not
Freudenthal, but the applied mathematician G. Braun, who died one year
before Hans' thesis work was finished and about whom I have been able
to find little information.

The subject of our conference, Poisson geometry, is now firmly
ensconced on both sides of the Atlantic as well as on both sides of
the equator, as Alan Weinstein observed this week, and Henrique
Bursztyn has kindly asked me to speak about Hans Duistermaat's impact
on the field.  In a narrow sense Duistermaat contributed very little
to Poisson geometry.  The subject dearest to his heart was
differential equations, although he had an unusual geometric
intuition.  As far as I know (thanks to Rui Loja Fernandes), the
notion of a \emph{Poisson manifold} appears just once in his written
work, namely in a book \cite{duistermaat;discrete-integrable} on
discrete dynamical systems on elliptic surfaces, which he finished not
long before his death and which has just been published.  But Poisson
\emph{brackets} can be found in most of his papers, and the fact is
that he has contributed many original ideas to the area.

\subsubsection*{Bispectral problem for Schr\"odinger equations}

For instance, his paper
\cite{duistermaat-grunbaum;differential-spectral} with Alberto
Gr\"unbaum continues to be influential in the literature on integrable
systems and noncommutative algebraic geometry.  It contains the
solution of the at first rather strange-sounding bispectral problem,
for what potentials $V(x)$ do the solutions $f(x,\lambda)$ of the
equation $-f''+Vf=\lambda f$ satisfy a differential equation in the
spectral parameter $\lambda$?  A comment on the bibliography of this
paper: most mathematicians disregard the classics, but Hans was never
afraid to go back to the sources.  He was widely read in the older
literature on analysis and differential geometry, and used it to great
effect in his own writings.

\subsubsection*{Resonances}

To mention a lesser-cited paper \cite{duistermaat;nonintegrability},
let me remind you of Duistermaat's insight how resonances in a
Hamiltonian system may preclude complete integrability, as explained
to us earlier this week by Nguyen Tien Zung.

\subsubsection*{Lie III}

Sometimes the flow of ideas was remarkably indirect.  His and Joop
Kolk's book on Lie groups \cite{duistermaat-kolk;lie-groups;;2000}
(which, despite having appeared in Springer's Universitext series, is
not exactly an elementary graduate text) is much closer to Lie's
original point of view pertaining to differential equations than
modern treatments such as Bourbaki \cite{bourbaki;groupes-algebres},
which are more algebraic in spirit.  Nevertheless the book is notable
for several innovations, particularly its proof of Lie's third
fundamental theorem in global form, which I think deserves to become
the standard argument and which runs in outline as follows.

The thing to be proved is that for every finite-dimensional real Lie
algebra $\g$ there exists a simply connected finite-dimensional real
Lie group whose Lie algebra is isomorphic to $\g$.  Choose a norm on
$\g$.  Then the space $P(\g)$ of continuous paths $[0,1]\to\g$
equipped with the supremum norm is a Banach space.  For each $\gamma$
in $P(\g)$ let $A_\gamma$ be the continuously differentiable path of
linear endomorphisms of $\g$ determined by the linear initial-value
problem
$$
A_\gamma'(t)=\ad(\gamma(t))\circ A_\gamma(t),\qquad
A_\gamma(0)=\id_\g.
$$
Duistermaat and Kolk prove that the multiplication law
$$(\gamma\cdot\delta)(t)=\gamma(t)+A_\gamma(t)\delta(t)$$
turns $P(\g)$ into a Banach Lie group.  Next they introduce a subset
$P(\g)_0$, which consists of all paths $\gamma$ that can be connected
to the constant path $0$ by a family of paths $\gamma_s$ which is
continuously differentiable with respect to $s$ and has the property
that
$$
\int_0^1A_{\gamma_s}(t)^{-1}\frac{\partial\gamma_s}{\partial
s}(t)\,dt=0
$$
for $0\le s\le1$.  The subset $P(\g)_0$ is a closed connected normal
Banach Lie subgroup of $P(\g)$ of finite codimension, and the quotient
$P(\g)/P(\g)_0$ is a simply connected Lie group with Lie algebra $\g$!

One of the many virtues of this proof is that it is manifestly
functorial: a Lie algebra homomorphism $\g\to\h$ induces a continuous
linear map on the path spaces $P(\g)\to P(\h)$, which is a
homomorphism of Banach Lie groups and maps the subgroup $P(\g)_0$ to
$P(\h)_0$, and therefore descends to a Lie group homomorphism.

For the details I refer you to the book; also be sure to read the
historical and bibliographical notes at the end of the chapter.  The
global form of Lie's third theorem appears to be due to \'E. Cartan
\cite{cartan;troisieme-theoreme-fondamental}, whose first proof was
based on the Levi decomposition.  A later version
\cite{cartan;topologie-groupes-lie} goes approximately as follows: to
build a simply connected group $G$ corresponding to $\g$, start with
the universal covering group $G_0$ of the adjoint group of $\g$, and
then construct $G$ as a central extension of $G_0$ by the centre of
$\g$ (viewed as an abelian Lie group).  The extension is obtained from
a cocycle $G_0\times G_0\to\lie{z}(\g)$, which Cartan finds by
integrating the infinitesimal cocycle that corresponds to the Lie
algebra extension $\lie{z}(\g)\to\g\to\g_0$.  This argument has a
natural interpretation in the language of differentiable group
cohomology, as was shown by Van Est
\cite{vanest;group-cohomology-lie-algebra}.  Much earlier Lie himself
suggested a different proof: he surmised that the Lie algebra $\g$
ought to be linear and that a group with Lie algebra $\g$ can
therefore be realized as a subgroup of an appropriate general linear
group.  This line of argument was justified a few years after Cartan
by Ado \cite{ado;representation-lie-algebras}.

The point of this for Poisson geometry is that a few years after
publication the Duistermaat-Kolk proof became at Alan Weinstein's
suggestion a central feature of Marius Crainic and Rui Loja Fernandes'
resolution of two longstanding problems in differential geometry: the
integrability of a Lie algebroid to a Lie groupoid
\cite{crainic-fernandes;integrability-lie}, and the integrability of a
Poisson manifold to a symplectic groupoid
\cite{crainic-fernandes;integrability-poisson}.  Curiously, a very
different work of Duistermaat, which I will get to later, also
impinges on these integrability problems.  My Cornell colleague
Leonard Gross has a paper in preparation that adapts the
Duistermaat-Kolk argument to certain infinite-dimensional situations.

Let me now discuss in a bit more detail four of Hans Duistermaat's
papers that are of obvious relevance to the topics of this conference,
namely those on the spectrum of elliptic operators
\cite{duistermaat-guillemin;spectrum-bicharacteristics}, global
action-angle variables \cite{duistermaat;global-action-angle}, the
quantum-mechanical spherical pendulum
\cite{duistermaat-cushman;quantum-spherical}, and the
Duistermaat-Heckman theorem \cite{duistermaat-heckman;variation}.

\section{The spectrum of positive elliptic operators and periodic 
bicharacteristics}

This paper \cite{duistermaat-guillemin;spectrum-bicharacteristics},
coauthored with Victor Guillemin, is Hans Duistermaat's most cited
work according to MathSciNet.  It is perhaps also his technically most
accomplished paper.  The authors consider a compact $n$-dimensional
manifold $X$ and a scalar elliptic pseudodifferential operator $P$ of
order $1$ on $X$ which is positive selfadjoint.  The spectrum of this
operator is a discrete set
$$
0\le\lambda_0\le\lambda_1\le\cdots\le\lambda_j\le\cdots\longto\infty.
$$
By placing a Dirac measure at each eigenvalue we obtain the spectral
distribution $\sigma_P=\sum_{j=0}^\infty\delta_{\lambda_j}$.  The
principal symbol $p$ of $P$ defines a Hamiltonian vector field $H_p$
on the punctured cotangent bundle $T^*X\setminus X$.

\begin{example*}
The main example to keep in mind is that of the Laplacian $\Delta$
defined with respect to a Riemannian metric on $X$.  For a suitable
constant $c$ the operator $c-\Delta$ is positive selfadjoint, so the
spectral theorem enables us to define a positive square root
$P=\sqrt{c-\Delta}$, which is pseudodifferential of order $1$ and
whose principal symbol is given by $p(x,\xi)=\norm{\xi}$ for $x\in X$
and $\xi\in T_xX$.  The Hamiltonian flow of $p$ is the geodesic spray
of $X$ (at least on the unit sphere bundle).
\end{example*}

The operator $P$ is a quantization of the classical observable $p$.
As explained for example in
\cite{guillemin;lectures-spectral-elliptic}, the classical analogue of
an eigenfunction of $P$ is a periodic trajectory of the Hamiltonian
vector field $H_p$ and the classical analogue of the eigenvalue is the
energy of the trajectory.  The periods of $H_p$ form a ``lattice''
which is dual to the set of eigenvalues of $P$.  The purpose of the
paper is to make this analogy precise.

Here I must interrupt myself to state that Hans never spoke to me in
such terms.  Many of us conceive of mathematics as a system of
grandiose functorial schemes and profound analogies or correspondences
suggested by the mysteries of nature.  Hans' mind worked differently
and the word ``quantization'' never crossed his lips except in jest.
Once he told me that on a visit to Moscow early in his career Gelfand
asked him what were his chief mathematical goals in life, and he had
no idea what to say.

What moved Hans Duistermaat, as far as I can see, was a gregarious and
competitive spirit that took him from one collaboration to the next
and from one mathematical problem to the next.  These are the same
qualities that made him a keen chess player, strong enough to have
once played former world champion A.~Karpov to a draw in a
simultaneous match.  An early manifestation of this spirit was his
eager participation in the sport of kite flying during his boyhood in
the Dutch East Indies, now the Republic of Indonesia.  The local
variant of the entertainment, which Wikipedia tells me is known as
kite fighting, required coating the flying line with glass and
abrasives for the purpose of ruining one's playmates' equipment.
Victor Guillemin relates in his acceptance notice for the 2003 Steele
Prize how Hans warned him, not for nothing, against getting involved
with a \emph{duistere maat} (murky companion).

Let us turn back to microlocal analysis and look at
$$
\hat{\sigma}_P(t)=\ca{F}\sigma_P=\sum_{j=0}^\infty e^{-i\lambda_jt},
$$
the Fourier transform of the spectral distribution $\sigma_P$.  This
can be seen as the distributional trace of the unitary operator
$e^{-itP}$, which is a Fourier integral operator.  A preliminary
result says that $\hat{\sigma}_P$ is a tempered distribution, and
therefore so is $\sigma_P$.  In particular the eigenvalue counting
function
$$
N_P(\lambda)=\sharp\{\,j\mid\lambda_j\le\lambda\,\}
$$
does not grow faster than a power of $\lambda$, which foreshadows
Weyl's law.  The next result is a first hint at the connection between
periods and eigenvalues.

\begin{theorem*}
$\hat{\sigma}_P$ is $C^\infty$ outside the set of periods of periodic
trajectories of $H_p$.
\end{theorem*}

In other words, the singular support of $\hat{\sigma}_P$ is contained
in the set of periods of the Hamiltonian $p$.  This is suggestive of
the Poisson summation formula, where one Fourier transforms a sum of
delta functions supported on a lattice and finds a sum of delta
functions supported on the dual lattice.  The results of Duistermaat
and Guillemin describe the singularities of $\hat{\sigma}_P$ and can
be viewed as a generalization of this elementary fact.

Every orbit is periodic of period $0$, so one expects $\hat{\sigma}_P$
to have a big singularity at $t=0$.  To focus on this singularity take
a smooth function $\chi$ such that $\hat{\chi}=\ca{F}\chi$ is a bump
function equal to $1$ in a small neighbourhood of $0$.  Then
$$
\hat{\chi}(t)\hat{\sigma}_P(t)=\sum_{j=0}^\infty
e^{-i\lambda_jt}\hat{\chi}(t)
=\ca{F}\biggl(\sum_j\chi(\lambda-\lambda_j)\biggr).
$$

\begin{theorem*}
We have an asymptotic expansion
$$
\sum_j\chi(\lambda-\lambda_j)\sim\frac1{(2\pi)^n}\sum_{k=0}^\infty
c_k\lambda^{n-1-k}
$$
as $\lambda\to\infty$.  The constants $c_k$ are independent of $\chi$.
The leading coefficient is
$$
c_0=\vol\{(x,\xi)\in T^*X\mid p(x,\xi)=1\}.
$$
\end{theorem*}

This yields all sorts of information about the spacing of the
eigenvalues, for example the following version of Weyl's law, which
says that the volume of phase space is asymptotically proportional to
the number of eigenvalues.

\begin{theorem*}
We have an asymptotic expansion
$$
N_P(\lambda)=\frac{a}{(2\pi)^n}\lambda^n+O(\lambda^{n-1})
$$
as $\lambda\to\infty$, where $a=\vol\{(x,\xi)\mid p(x,\xi)\le1\}$.
\end{theorem*}

A further analysis leads to a ``residue formula'', which describes the
poles of $\hat{\sigma}_P$ at nonzero periods.

\begin{theorem*}
Let $T\ne0$.  Assume that all periodic orbits of $H_p$ of period $T$
are isolated and nondegenerate.  Then
$$
\lim_{t\to
T}(t-T)\hat{\sigma}_P(t)=\sum_\gamma\frac{T_{0,\gamma}}{2\pi}
\frac{i^{m_\gamma}}{\bigabs{\det(I-d\Pi_\gamma)}^{1/2}}.
$$
\end{theorem*}

There is a close resemblance between this formula and the Lefschetz
formula for elliptic complexes of Atiyah and Bott
\cite{atiyah-bott;lefschetz-fixed-elliptic-complex}.  The sum on the
right is over all closed orbits $\gamma$ of period $T$; $T_{0,\gamma}$
is the primitive period of $\gamma$; and $\Pi_\gamma$ is the
Poincar\'e return map of $\gamma$.  Nondegeneracy of $\gamma$ means
that $\det(I-d\Pi_\gamma)\ne0$.  The integer $m_\gamma$ is a Maslov
index.  For $P=\sqrt{c-\Delta}$ it is the Morse index of the geodesic
$\gamma$ for the Euler-Lagrange functional.  Eckhard Meinrenken showed
in an early paper \cite{meinrenken;conley-zehnder} that for general
$P$ the number $m_\gamma$ can be interpreted as a Conley-Zehnder
index.

Duistermaat and Guillemin did not do this work in isolation.  Some of
the most important prior mathematical work on the subject is that of
Weyl \cite{weyl;verteilungsgesetz}, which was inspired by Planck's
model of black-body radiation, and H\"ormander
\cite{hormander;spectral-elliptic;acta}.  Roughly contemporaneous work
includes that of Gutzwiller \cite{gutzwiller;periodic-quantization},
Colin de Verdi\`ere \cite{colin-de-verdiere;spectre-laplacien}, and
Chazarain \cite{chazarain;formule-de-poisson}.  See \cite{arendt;weyl}
for a historical survey that takes in a good deal of the physics
literature.  For later developments the reader can consult the Fourier
volume in honour of Colin de Verdi\`ere, particularly Colin's own
contribution \cite{colin-de-verdiere;spectrum-laplace} to that volume.

\section{On global action-angle variables}

I will give a slightly anachronistic account of Duistermaat's paper
\cite{duistermaat;global-action-angle} on monodromy in integrable
systems, which takes into consideration later work of Dazord and
Delzant \cite{dazord-delzant;actions-angles}.

Let $B$ be a connected $n$-manifold.  A \emph{Lagrangian fibre bundle}
over $B$ is a triple $\ca{L}=(M,\omega,\pi)$, where $(M,\omega)$ is a
symplectic $2n$-manifold and $\pi\colon M\to B$ a surjective
submersion with Lagrangian fibres.  To keep things simple we will
assume that the fibres of $\pi$ are compact and connected.

The standard way to obtain such a bundle is to start with an
integrable Hamiltonian system and throw out the singularities of the
energy-momentum map.  The simplest Lagrangian fibre bundle over a
given base $B$ is as follows.

\begin{example*}
Let $p\colon B\to\R^n$ be a local diffeomorphism.  Let $\T$ be the
circle $\R/\Z$.  The \emph{angle form} on $\T$ is $dq$, where $q$ is
the coordinate on $\R$.  Let $M$ be the product $B\times\T^n$,
equipped with the symplectic form $\omega=\sum_{j=1}^ndp_j\wedge
dq_j$.  Let $\pi\colon M\to B$ be the projection onto the first
factor.  The functions $p_j$ are the \emph{action variables} and the
(multivalued) functions $q_j$ are the \emph{angle variables}.  The map
$p\circ\pi$ is a momentum map for the translation action of $\T^n$ on
the second factor of $M$.
\end{example*}

An \emph{isomorphism} of Lagrangian fibre bundles over $B$ is given by
a symplectomorphism of the total spaces that induces the identity map
on the base.  There is an equally obvious notion of
\emph{localization}, that is restriction of a Lagrangian fibre bundle
to an open subset of the base.  The Liouville-Mineur-Arnold theorem
states that every Lagrangian fibre bundle admits local action-angle
variables, i.e.\ is locally isomorphic to $B\times\T^n$.  The problem
solved by Duistermaat, and before him in special cases by
Nehoro{\v{s}}ev \cite{nehorosev;action-angle}, is when a Lagrangian
fibre bundle over $B$ admits global action-angle variables, i.e.\ is
globally isomorphic to $B\times\T^n$.

\begin{theorem*}
A Lagrangian fibre bundle admits global action-angle variables if and
only if two invariants, $\mu(P)$ (the affine monodromy) and
$\lambda(\ca{L})$ (the Lagrangian class), vanish.
\end{theorem*}

(This is not quite the formulation given by Duistermaat.  The
Lagrangian class was introduced later by Dazord and Delzant.  The
quantity called monodromy by Duistermaat is what I will call here the
linear monodromy, which ignores the translational part of the affine
monodromy.  See below for a full discussion.)

Just as interesting as this theorem is the fact that many commonplace
integrable systems do \emph{not} admit global action-angle variables,
for instance Huygens' spherical pendulum, which Duistermaat analyses
in detail.  Let me now explain the two invariants.

\subsection*{Monodromy}

Let $\ca{L}=(M,\omega,\pi)$ be a Lagrangian fibre bundle over $B$.
The map $TM\to T^*M$ given by $v\mapsto\iota(v)\omega$ is a bundle
isomorphism, and we denote its inverse by $\omega^\sharp\colon T^*M\to
TM$.  Let $m\in M$ and put $b=\pi(m)\in B$.  Given a covector
$\alpha\in T_b^*B$, the projection and the symplectic form produce a
tangent vector $v_m(\alpha)$,
$$
T_b^*B\overset{\pi^*}\longto T_m^*M\overset{\omega^\sharp}\longto
T_mM,\qquad\alpha\longmapsto\pi^*(\alpha)\longmapsto
\omega^\sharp\pi^*(\alpha)=v_m(\alpha).
$$
Since we can write $\alpha=d_bf$ for a suitable function $f$, we see
that $v_m(\alpha)$ is the value at $m$ of the Hamiltonian vector field
$H_{\pi^*f}$, and therefore is tangent to the fibre $\pi^{-1}(b)$.
The fibre being compact, the vector field $v(\alpha)$ is complete, and
we denote by $\phi_b(\alpha)\colon\pi^{-1}(b)\to\pi^{-1}(b)$ its time
$1$ flow.  The map
$$\phi_b\colon T_b^*B\times\pi^{-1}(b)\longto\pi^{-1}(b)$$
defined by $\phi_b(\alpha,m)=\phi_b(\alpha)(m)$ is an action of the
abelian Lie group $T_b^*B\cong\R^n$ on $\pi^{-1}(b)$.  The map
$\alpha\mapsto v_m(\alpha)$ is an isomorphism $T_b^*B\to
T_m(\pi^{-1}(b))$, so, the fibre $\pi^{-1}(b)$ being connected, we
conclude that the action $\phi_b$ is transitive and locally free.  The
kernel of the action $P_b\cong\Z^n$ is the \emph{period lattice}
at $b$.  Collecting these fibrewise actions gives us an action
$$\phi\colon T^*B\times_BM\longto M$$
of the bundle of Lie groups $T^*B\to B$ on the bundle $M\to B$.  The
kernel of this bundle action is the bundle of free abelian groups
$P=\coprod_bP_b$ over $B$, called the \emph{period bundle}.  The
fibrewise quotient
$$T=T^*B/P$$
is a bundle over $B$ with general fibre the torus $\T^n$ and structure
group $\Aut(\T^n)\cong\GL(n,\Z)$.  The quotient action
$$\phi_T\colon T\times_BM\longto M,$$
which we will write as $\phi_T(t,m)=t\cdot m$, makes $M$ a
\emph{$T$-torsor}, a principal homogeneous space for the torus bundle
$T$ in the sense that the map $T\times_BM\to M\times_BM$ defined by
$(t,m)\mapsto(m,t\cdot m)$ is a diffeomorphism.

The $T^*B$-action defines for each $1$-form on the base
$\alpha\in\Omega^1(B)$ a diffeomorphism $\phi(\alpha)$ from $M$ to
itself which induces the identity on $B$.  This diffeomorphism
transforms the symplectic form as follows.

\begin{lemma*}
$\phi(\alpha)^*\omega=\omega+\pi^*d\alpha$ for every
$\alpha\in\Omega^1(B)$.
\end{lemma*}

Recall that a Lagrangian section of a cotangent bundle is the same as
a closed $1$-form.  Since sections of $P$ induce the identity map
on $M$, the lemma tells us therefore that $P$ is a Lagrangian
submanifold of $T^*B$.

This has various desirable consequences.  First of all, applying the
lemma to the translation action of $T^*B$ on itself we conclude that
the standard symplectic form is preserved by the $P$-action and so
descends to a symplectic form $\omega_T$ on $T$.  Thus the Lie group
bundle $T$ itself is a Lagrangian fibre bundle over $B$.

More importantly, we see that on any sufficiently small open subset
$U$ of the base there exists a coordinate system
$p=(p_1,p_2,\dots,p_n)$ such that
$$\ca{F}(p)=(dp_1,dp_2,\dots,dp_n)$$
is a frame of the local system $P|U$.  These preferred coordinate
systems determine what following recent usage I will call a
\emph{tropical affine structure} on $B$, that is an atlas with values
in the pseudogroup defined by the \emph{tropical affine group}
$\G=\GL(n,\Z)\ltimes\R^n$.  (See e.g.\
\cite[Chapter~1]{gross;tropical-mirror}.)  Conversely, this atlas
determines the Lagrangian lattice bundle $P$.

Analytic continuation of the coordinate system $p$ along a loop
$\gamma$ in $B$ based at $b\in U$ gives a new coordinate system $p'$
at $b$, which is related to $p$ by a transformation $g_\gamma\in\G$.
The corresponding local frames $\ca{F}(p)$ and $\ca{F}(p')$ of
$P$ are related by the linear part
$g_{0,\gamma}\in\G_0=\GL(n,\Z)$ of the affine transformation
$g_\gamma$.  The map $\gamma\mapsto g_\gamma$ induces a homomorphism
from $\pi_1(B,b)$ to $\G$.  The conjugacy class of this homomorphism,
$$
\mu(P)\in\Hom(\pi_1(B),\G)/\Ad(\G)\cong H^1(B,\G),
$$
is the \emph{affine monodromy} of $P$.  (Here $H^1(B,\G)$ denotes the
cohomology set of $B$ with coefficients in the group $\G$ equipped
with the discrete topology.)  The conjugacy class defined by the map
$\gamma\mapsto g_{0,\gamma}$,
$$
\mu_0(P)\in\Hom(\pi_1(B),\G_0)/\Ad(\G_0)\cong H^1(B,\G_0),
$$
is the \emph{linear monodromy}, which determines the isomorphism class
of the local system $P$.  The monodromy depends only on the
affine structure of $B$, not on $M$ or its symplectic structure.

The linear monodromy $\mu_0(P)$ is trivial if and only if the local
system $P$ is trivial.  In that case $T\cong B\times\T^n$ is
isomorphic to a trivial bundle of Lie groups, $M$ is a principal
$\T^n$-bundle over $B$, and $P$ has a global frame of closed $1$-forms
$(\alpha_1,\alpha_2,\dots,\alpha_n)$.  We can then find a covering
$f\colon\tilde{B}\to B$ of the base and a local diffeomorphism
$\tilde{p}\colon\tilde{B}\to\R^n$ such that
$f^*\alpha_j=d\tilde{p}_j$.  The full monodromy $\mu(P)$ is trivial if
and only if $P$ is trivial and the $\alpha_j$ are exact.  If that
happens we can define global single-valued action variables $p\colon
B\to\R^n$, and $p\circ\pi$ is a momentum map for the $\T^n$-action on
$M$.

\subsection*{Chern class and Lagrangian class}

The existence of global angle variables on a Lagrangian fibre bundle
$\ca{L}=(M,\omega,\pi)$ is tantamount to the existence of a global
Lagrangian section of $\pi\colon M\to B$.

First let us consider plain smooth sections of $\pi$.  We need to
introduce a few sheaves of abelian groups on the base space $B$.
There is $\Omega^k$, the sheaf of smooth $k$-forms, and its subsheaf
$\ca{Z}^k$ of closed $k$-forms.  Then there is the sheaf of smooth
sections of $T$, which we will call $\ca{T}$, and the sheaf of locally
constant sections of $P$, which we will call $\ca{P}$.  Let
$\{U_i\}_{i\in I}$ be an open cover of $B$ and suppose that we have
local smooth sections $s_i\colon U_i\to M$ of $\pi$.  Since $M$ is a
$T$-torsor, over each intersection $U_{ij}=U_i\cap U_j$ we have a
unique section $t_{ij}\in\ca{T}(U_{ij})$ such that $s_i=t_{ij}\cdot
s_j$.  The tuple $t=(t_{ij})$ is a \v{C}ech $1$-cocycle and defines an
element $[t]\in H^1(B,\ca{T})$.

Since $T$ is the quotient bundle $T^*B/P$, on the level of
sheaves we have a short exact sequence
$$
0\longto\ca{P}\longto\Omega^1\longto\ca{T}\longto0.
$$
The sheaf $\Omega^1$ is fine, so the long exact cohomology sequence
gives canonical isomorphisms
$$H^k(B,\ca{T})\cong H^{k+1}(B,\ca{P})$$
for all $k\ge0$.  The image $c(M)\in H^2(B,\ca{P})$ of $[t]$ is the
\emph{Chern class} of the $T$-torsor $M$ and it is the
obstruction to the existence of a global section of $\pi$.  It is
independent of the symplectic structure on $M$.

Since $P$ is Lagrangian, the sheaf $\ca{P}$ is a subsheaf of
$\ca{Z}^1$, and therefore the exterior derivative
$d\colon\Omega^1\to\ca{Z}^2$ descends to a morphism
$$d_P\colon\ca{T}\longto\ca{Z}^2.$$
A section $t$ of $T$ is \emph{closed} if $d_Pt=0$.  If the open sets
$U_i$ are small enough, we can choose the local sections $s_i$ to be
Lagrangian, which implies that the transition functions $t_{ij}$ are
closed.  Thus the $t_{ij}$ are sections of the subsheaf
$\ca{K}=\ker(d_P)$ of $\ca{T}$, and the corresponding cohomology class
lives in $H^1(B,\ca{K})$.  This is the \emph{Lagrangian class}
$\lambda(\ca{L})$, which is implicit in the paper of Dazord and
Delzant but was named by Zung
\cite{zung;topology-integrable-hamiltonian}, and it is the obstruction
to the existence of a global \emph{Lagrangian} section of $\pi$.
Given a Lagrangian section $s$, the map $T^*B\to M$ defined by
$(b,\alpha)\mapsto\phi(\alpha)(s(b))$ identifies the Lagrangian fibre
bundle $T$ with $\ca{L}$.

Therefore the vanishing of the Lagrangian class $\lambda(\ca{L})$ is
equivalent to $\ca{L}$ being isomorphic as a Lagrangian fibre bundle
to $T$.  The vanishing of both $\lambda(\ca{L})$ and the affine
monodromy $\mu(P)$ is equivalent to the existence of global
action-angle variables.  This is the version of Duistermaat's theorem
established by Dazord and Delzant (who, by the way, also considered
the case of less than fully integrable systems).

\subsection*{Symplectic torsors}

Dazord and Delzant went on to prove that the Lagrangian class
completely classify all Lagrangian fibre bundles on a tropical affine
manifold.  Let us widen our view a little by fixing a tropical affine
manifold $B$ with period bundle $P$ and torus bundle $T=T^*B/P$, and
examining arbitrary $T$-torsors over $B$.  Any such torsor $\pi\colon
M\to B$ has a well-defined Chern class $c(M)\in H^2(B,\ca{P})$, where
as before $\ca{P}$ is the locally constant sheaf of sections of $P$.
In fact, just as for principal bundles the cohomology group
$H^2(B,\ca{P})$ classifies $T$-torsors up to isomorphism.  (If the
linear monodromy $\mu_0(P)$ vanishes, then $P$ is the constant local
system $\Z^n$, a $T$-torsor is an ordinary principal $\T^n$-bundle,
and the Chern class is the ordinary Chern class in $H^2(B,\Z^n)$.)

Let us think about all possible symplectic forms $\omega$ on $M$ which
vanish on the fibres of $\pi$, so making $\ca{L}=(M,\omega,\pi)$ into
a Lagrangian fibre bundle.  We will call $\ca{L}$ a \emph{symplectic
$T$-torsor} with \emph{total space} $M$.  As before we regard two
symplectic $T$-torsors as isomorphic if the total spaces are
symplectomorphic via a diffeomorphism that fixes the base $B$.  The
collection of isomorphism classes $[\ca{L}]$ is an analogue of the
Picard group of an algebraic variety and we will denote it by
$\Pic(B,P)$.

The set $\Pic(B,P)$ is equipped with two algebraic operations.  The
\emph{opposite} of $\ca{L}=(M,\omega,\pi)$ is
$-\ca{L}=(M,-\omega,\pi)$.  (Negating the symplectic form has the
effect of reversing the $T$-action, i.e.\ composing it with the
automorphism $t\mapsto t^{-1}$ of $T$.)  Given two symplectic
$T$-torsors $\ca{L}_1=(M_1,\omega_1,\pi_1)$ and
$\ca{L}_2=(M_2,\omega_2,\pi_2)$, define $M$ to be the $T$-torsor
$(M_1\times_BM_2)/T^-$, where $T^-$ is the antidiagonal subbundle
$\{(t,t^{-1})\mid t\in T\}$ of $T\times_BT$.  It is a theorem of Ping
Xu \cite{xu;morita-groupoid} that the form $\omega_1+\omega_2$ on
$M_1\times_BM_2$ descends to a symplectic form $\omega$ on $M$ which
makes $\ca{L}=(M,\omega,\pi)$ into a symplectic $T$-torsor.  We call
$\ca{L}$ the \emph{sum} of $\ca{L}_1$ and $\ca{L}_2$.  The operation
$[\ca{L}_1]+[\ca{L}_2]=[\ca{L}_1+\ca{L}_2]$ turns $\Pic(B,P)$ into an
abelian group.  The zero element is $[T]$ and the opposite of
$[\ca{L}]$ is $[-\ca{L}]$.

Can we explicitly describe the Picard group $\Pic(B,P)$?  The
Poincar\'e lemma implies that
$$
0\longto\ca{Z}^k\longto\Omega^k\overset{d}\longto\ca{Z}^{k+1}\longto0
$$
is a short exact sequence of sheaves.  To begin with, this gives us
isomorphisms
$$H^l(B,\ca{Z}^k)\cong H^{k+l}(B,\R),$$
because $\Omega^k$ is fine.  Furthermore, taking $k=1$ and dividing
the first two terms by $\ca{P}$ we get the short exact sequence
$$
0\longto\ca{Z}^1/\ca{P}\longto\ca{T}\overset{d_P}\longto
\ca{Z}^2\longto0.
$$
This identifies the kernel $\ca{K}=\ker(d_P)$ with $\ca{Z}^1/\ca{P}$
and yields a long exact sequence
\begin{multline*}
0\longto H^0(B,\ca{K})\longto H^0(B,\ca{T})\overset{d_{P,*}}\longto
H^0(B,\ca{Z}^2)\overset{\partial}\longto H^1(B,\ca{K})\longto
H^1(B,\ca{T})
\\
\overset{d_{P,*}}\longto H^1(B,\ca{Z}^2)\overset{\partial} \longto
H^2(B,\ca{K})\longto\cdots
\end{multline*}
Substituting $H^k(B,\ca{Z}^2)\cong H^{k+2}(B,\R)$ and
$H^k(B,\ca{T})\cong H^{k+1}(B,\ca{P})$, and noticing that
$H^k(B,\ca{K})\to H^k(B,\ca{T})\cong H^{k+1}(B,\ca{P})$ is the
connecting homomorphism $\delta$ for the short exact sequence
$$
0\longto\ca{P}\longto\ca{Z}^1\longto\ca{K}\longto0,
$$
we obtain the long exact sequence that we want,
\begin{multline*}
0\longto H^0(B,\ca{K})\overset{\delta}\longto
H^1(B,\ca{P})\overset{d_{P,*}}\longto
H^2(B,\R)\overset{\partial}\longto
H^1(B,\ca{K})\overset{\delta}\longto H^2(B,\ca{P})
\\
\overset{d_{P,*}}\longto H^3(B,\R)\overset{\partial}\longto
H^2(B,\ca{K})\overset{\delta}\longto\cdots
\end{multline*}

If a $T$-torsor $M$ admits a symplectic form $\omega$ vanishing on the
fibres, then $\delta$ maps the Lagrangian class
$\lambda(M,\omega,\pi)$ to the Chern class $c(M)$, and therefore
$d_{P,*}c(M)=0$.  So we see that $d_{P,*}c(M)=0$ is a necessary
condition for $M$ to be the total space of a symplectic $T$-torsor.
Dazord and Delzant show that this condition is actually sufficient,
and that every $\lambda\in H^1(B,\ca{K})$ satisfying
$\delta\lambda=c(M)$ is the Lagrangian class of a unique isomorphism
class of Lagrangian fibre bundles with total space $M$.  The
conclusion is as follows.

\begin{theorem*}
Let $B$ be a tropical affine manifold with period bundle $P$.  Let $T$
be the torus bundle $T^*B/P$, let $\ca{T}$ be the sheaf of smooth
sections of $T$, and let $\ca{K}$ be the kernel of the sheaf
homomorphism $d_P\colon\ca{T}\to\ca{Z}^2$.
\begin{enumerate}
\item
The map $\Pic(B,P)\to H^1(B,\ca{K})$ defined by
$[\ca{L}]\mapsto\lambda(\ca{L})$ is a group isomorphism.
\item
We have a short exact sequence
$$
0\longto H^2(B,\R)/d_{P,*}H^1(B,\ca{P})\longto\Pic(B,P)\longto
\ker(d_{P,*})\longto0.
$$
\end{enumerate}
\end{theorem*}

This theorem gives us two different descriptions of the identity
component of the Picard group, namely $\Pic^0(B,P)$ is equal to
the group of symplectic torsors of ``degree'' (i.e.\ Chern class) $0$,
and
$$\Pic^0(B,P)\cong H^2(B,\R)/d_{P,*}H^1(B,\ca{P}).$$
The ``N\'eron-Severi group'' (i.e.\ component group)
$\Pic(B,P)/\Pic^0(B,P)$ is isomorphic to the subgroup $\ker(d_{P,*})$
of $H^2(B,\ca{P})$.  If the base $B$ is of finite type, the Picard
group is finite-dimensional and the N\'eron-Severi group is finitely
generated.

Suppose that we are given a $T$-torsor $\pi\colon M\to B$ with Chern
class $c\in H^2(B,\ca{P})$ and let us denote by
$\Pic(M,B,P)\cong\delta^{-1}(c)$ the collection of isomorphism classes
of symplectic $T$-torsors with total space $M$.  The theorem tells us
that $\Pic(M,B,P)$ is nonempty if and only if $d_{P,*}c=0$ and that
the group $\Pic^0(B,P)$ acts simply transitively on $\Pic(M,B,P)$.
The $\Pic^0(B,P)$-action is the \emph{gauge action} given by the
formula $[\sigma]\cdot[M,\omega,\pi]=[M,\omega+\pi^*\sigma,\pi]$,
where $\sigma\in Z^2(B)$ is a de Rham representative of a class in
$H^2(B,\R)$.

Zung has obtained a version of these results for certain
\emph{singular} Lagrangian fibrations.

\subsection*{Twisted symplectic torsors}

It is instructive to go one step further in the long exact sequence
and ask what happens if $d_{P,*}c(M)$ is nonzero.  This leads to a
``nonholonomic'' version of the Duistermaat-Dazord-Delzant theorems.
I will outline the results and publish the proofs elsewhere.  We
define a \emph{twisted Lagrangian fibre bundle}
$\ca{L}=(M,\omega,\pi)$ over a base manifold $B$ in the same way as a
Lagrangian fibre bundle, except that we relax the requirement that
$\omega$ be closed to the requirement that $d\omega$ be \emph{basic}
in the sense that $d\omega=\pi^*\eta$ for a $3$-form $\eta$ on $B$.
Thus $\omega$ is an almost symplectic form.  The form
$\eta=\eta(\ca{L})$ is closed and uniquely determined by $\omega$.  It
is an isomorphism invariant and we refer to it as the \emph{twisting
form} of the twisted Lagrangian fibre bundle $\ca{L}$.  The pair
$(\omega,\eta)$ is a cocycle in the relative de Rham complex of the
projection $\pi$.

It turns out that, just as in the Lagrangian case, the cotangent
bundle $T^*B$ acts on the total space $M$ of a twisted Lagrangian
fibre bundle $\ca{L}$ and that the kernel of the action is a bundle of
lattices $P$, which is Lagrangian with respect to the standard
symplectic form on $T^*B$.  So again $B$ is a tropical affine manifold
and $M$ is a torsor for the torus bundle $T=T^*B/P$.

We now fix the tropical affine manifold $(B,P)$ and look at any
twisted Lagrangian bundle $\ca{L}=(M,\omega,\pi)$ which is at the same
time a $T$-torsor.  We assume that the almost symplectic form $\omega$
is \emph{compatible} with the $T$-action in the sense that the
$T^*B$-action on $M$ induced by $\omega$ has kernel $P$.  We call such
an $\ca{L}$ a \emph{twisted symplectic $T$-torsor} and set ourselves
the task of classifying up to isomorphism all twisted symplectic
$T$-torsors.  We denote the set of isomorphism classes by
$\TPic(B,P)$.

The first observation is that this set is an abelian group in the same
way as the ordinary Picard group.  We will refer to $\TPic(B,P)$ as
the \emph{twisted} Picard group of the tropical affine manifold
$(B,P)$.

The next observation is that every $T$-torsor $M$ possesses a
compatible almost symplectic form and that the extent to which it is
not closed is measured by the class $d_{P,*}c(M)$.

\begin{theorem*}
Every $T$-torsor $\pi\colon M\to B$ possesses a compatible almost
symplectic form $\omega$.  Its twisting form $\eta\in Z^3(B)$
satisfies $[\eta]=d_{P,*}c(M)$.
\end{theorem*}

The Dazord-Delzant theorem generalizes as follows.

\begin{theorem*}
We have an exact sequence
$$
0\longto \Omega^2(B)/d_PH^0(B,\ca{T})\longto\TPic(B,P)\longto
H^2(B,\ca{P})\longto0.
$$
\end{theorem*}

So the moduli space of twisted Lagrangian fibre bundles is typically
infinite-dimensional.  These degrees of freedom can be taken away by
introducing a coarser form of gauge equivalence, namely by letting an
arbitrary $2$-form $\sigma\in\Omega^2(B)$ on the base act on a twisted
Lagrangian fibre bundle $\ca{L}=(M,\omega,\pi)$ by the formula
$\sigma\cdot\ca{L}=(M,\omega+\pi^*\sigma,\pi)$.  This action changes
the twisting form by the exact $3$-form $d\beta$.  It follows from the
theorem that $\TPic(B,P)/\Omega^2(B)$ is isomorphic to
$H^2(B,\ca{P})$, in other words every $T$-torsor has a compatible
almost symplectic structure which is unique up to coarse gauge
equivalence.

A twisted symplectic $T$-torsor does not have a well-defined
Lagrangian class, but the \emph{difference}
$\ca{L}_1-\ca{L}_2=\ca{L}_1+-\ca{L}_2$ of two twisted symplectic
$T$-torsors that have the same twisting forms,
$\eta(\ca{L}_1)=\eta(\ca{L}_2)$, is a symplectic $T$-torsor and
therefore has a well-defined Lagrangian class.  It follows that if we
fix a closed $3$-form $\eta\in Z^3(B)$ the set of isomorphism classes
of twisted symplectic $T$-torsors with twisting form $\eta$ is a
principal homogeneous space of $\Pic(B,P)$.  If in addition we fix a
class $c\in H^2(B,\ca{P})$ satisfying $d_{P,*}c=[\eta]$, then the set
of isomorphism classes of twisted symplectic $T$-torsors with Chern
class $c$ and twisting form $\eta$ is a principal homogeneous space of
$\Pic^0(B,P)$.

\subsection*{Groupoids and realizations}

At the end of my talk Alan Weinstein pointed out that Duistermaat's
study of global action-angle variables provided one of the incentives
for him to formulate the symplectic groupoid program
\cite{coste-dazord-weinstein;groupoides},
\cite{weinstein;symplectic-groupoids-poisson-manifolds}.  In the
language of that program a Lagrangian fibre bundle is nothing but a
\emph{realization} of the base manifold $B$ equipped with the zero
Poisson structure, and the torus bundle $T$ is a symplectic
groupoid over $B$ (with source and target maps being equal) which
\emph{integrates} this Poisson manifold.

Every manifold $B$ with zero Poisson structure is obviously integrable
and the associated source-simply connected symplectic groupoid is just
the cotangent bundle $T^*B$.  What makes the tropical affine case
special is the existence of a \emph{proper} symplectic groupoid $T$
which integrates the trivial Poisson structure.  The
Duistermaat-Dazord-Delzant theorems then amount to a classification of
all realizations of $B$ which are free and fibre-transitive under $T$.
The group $\Pic(B,P)$ is referred to as the \emph{static} Picard group
of the groupoid $T$ in \cite{bursztyn-weinstein;picard-poisson}.  (The
full, noncommutative, Picard group is the semidirect product of
$\Pic(B,P)$ with the group of tropical affine automorphisms of $B$.)
The twisted case also fits into this framework, as one can see from
the papers \cite{bursztyn-crainic-weinstein-zhu;integration-twisted}
and \cite{cattaneo-xu;integration-twisted}.

\subsection*{The quantum-mechanical spherical pendulum}

Having spent far more time on action-angle variables than I intended,
let me be very brief about the quantum-mechanical picture.  A
treatment of quantum monodromy in the spherical pendulum was given by
Richard Cushman and Hans Duistermaat
\cite{duistermaat-cushman;quantum-spherical}.  A different
interpretation was given soon afterwards by Victor Guillemin and
Alejandro Uribe \cite{guillemin-uribe;monodromy-quantum}.  Let me
quote from Hans' review of the latter paper in the \emph{Mathematical
Reviews}, which he starts by explaining his own approach:

\newenvironment{Quotation}{\begin{list}{}{%
\setlength\leftmargin{\parindent}\setlength\rightmargin{\leftmargin}}
\item[]\ignorespaces\small}{\unskip\end{list}}

\begin{Quotation}
If one considers the Schr\"odinger operator
$E=-(\hbar^2/2)\Delta+V$, where $\Delta$ is the Laplace operator on
the $2$-dimensional standard sphere $S$ in $\R^3$ and the potential
$V$ is the vertical coordinate function, then the rotational symmetry
around the vertical axis yields an operator
$L=i\hbar(x_1\partial/\partial x_2-x_2\partial/\partial x_1)$ which
commutes with $S$ [\emph{sic}].

Replacing $i\hbar\partial/\partial x_j$ by the conjugate variable
$p_j$, we get principal symbols $e$ and $l$ of $E$ and $L$,
respectively, which Poisson commute and define an integrable
Hamiltonian system on the phase space $T^*S$, the cotangent bundle of
$S$.  One has straightforward generalizations to $n$ commuting
operators $E_1$,\dots, $E_n$ with principal symbols $e_1$,\dots, $e_n$
on $n$-dimensional manifolds $M$.

Because the operators $E_j$ commute, one has common eigenfunctions
$\psi_k$, $k=1$, $2$,\dots, with eigenvalues $\eps_{j,k}$
($E_j\psi_k=\eps_{j,k}\cdot\psi_k$).  The rule for finding the
$n$-dimensional spectrum $(\eps_{1,k},\dots,\eps_{n,k})\in\R^n$ for
$k=1$, $2$,\dots, asymptotically for $\hbar\downarrow0$ and near a
regular value of the mapping $(e_1,\cdots,e_n)$ is as follows.  One
constructs locally so-called action variables, which are functions
$(a_1,\dots,a_n)$ of the $(e_1,\dots,e_n)$, in such a way that
$(\partial a_i/\partial e_j)$ is invertible and the Hamiltonian flows
of the $a_j$ are periodic with period $2\pi$.  Then the
$n$-dimensional spectrum is given asymptotically by
$a^{-1}(\Z^n+\alpha)$, where $\Z^n$ is the integer lattice, $\alpha$
is a Maslov shift, and $a$ is the vector of action variables given
above.  This means that the actions, in particular the nonexistence of
global action variables, can be read off from the asymptotics of the
spectrum.
\end{Quotation}

He proceeds to explain the different approach taken by Guillemin and
Uribe.  Although it has some very convincing illustrations, Cushman
and Duistermaat's paper is little more than an announcement and there
does not seem to exist a more comprehensive version.  Some ten years
after its appearance experimental evidence of quantum monodromy was
found and finally V\~u Ng\d{o}c San wrote two papers
\cite{vu;quantum-monodromy-integrable-systems},
\cite{vu;quantum-monodromy-bohr-sommerfeld} clarifying and elaborating
on Cushman and Duistermaat's ideas.

\section{Duistermaat-Heckman}

Of all Hans Duistermaat's accomplishments the best known to
differential geometers is probably the Duistermaat-Heckman theorem
\cite{duistermaat-heckman;variation}.  This is so familiar to most of
the audience that I passed it over in my talk, but in this written
version I can't resist making some remarks about it.  Recall that in
its simplest form the theorem states that
$$
\int_M\exp(\omega-tf)=\int_X\frac{\exp(\omega-tf)}{\e(X,t)}.
$$
Here $(M,\omega)$ is a compact symplectic manifold, $t$ is a complex
parameter, $f$ is a periodic Hamiltonian, $X$ is the critical manifold
of $f$, and $\e(X,t)$ is the equivariant Euler class of the normal
bundle of $X$ in $M$.  The integral on the left is to be interpreted
as the integral of $e^{-tf}\omega^n/n!$, where $2n=\dim M$.  This is
precisely the Fourier-Laplace transform of the measure
$f_*(\omega^n/n!)$ obtained by pushing forward the Liouville measure
$\omega^n/n!$ to the real line.  The critical manifold $X$ usually
consists of connected components of various dimensions, so the
integral on the right is to be read as a sum of integrals, one for
each component.

The theorem contains as a special case Archimedes' result that the
surface area of a sphere is equal to that of the circumscribed
cylinder, an illustration of which, according to Cicero \cite[Liber~V,
\S\S\,64--66]{cicero;disputations}, adorned the Syracusan's tomb.  A
modern antecedent of the theorem is Bott's residue formula for
holomorphic vector fields~\cite{bott;residue}.  Soon after publication
three interesting alternative proofs appeared, one based on the
localization principle in equivariant cohomology by Atiyah and Bott
\cite{atiyah-bott;moment-map-equivariant-cohomology} and Berline and
Vergne \cite{berline-vergne;zeros} (see also
\cite{berline-getzler-vergne;heat-kernels} and
\cite{guillemin-sternberg;supersymmetry-equivariant}), one based on
partial action-angle variables by Dazord and Delzant
\cite{dazord-delzant;actions-angles}, and one based on the coisotropic
embedding theorem by Guillemin and
Sternberg~\cite{guillemin-sternberg;birational}.

\subsection*{The index theorem}

My favourite meta-application of the Duistermaat-Heckman theorem is
Atiyah's heuristic derivation \cite{atiyah;circular-symmetry} of the
Atiyah-Singer index theorem for the Dirac operator suggested by ideas
of Witten \cite{witten;supersymmetry-morse}.  Let $X$ be a compact
Riemannian manifold and let $M=C^\infty(S^1,X)$ be the loop space of
$X$.  A tangent vector to $M$ at a loop $\gamma$ is a vector field
along $\gamma$, i.e.\ a section of $\gamma^*(TX)$.  The loop space has
a Riemannian structure: the inner product of two tangent vectors
$s_1$, $s_2\in T_\gamma M$ is defined to be the integral
$\int_{S^1}(s_1(\theta),s_2(\theta))\,d\theta$.  The circle $S^1$ acts
on $M$ by spinning the loops, and we let $\alpha$ be the $1$-form on
$M$ dual to the infinitesimal generator of this action.  Then
$\omega=d\alpha$ is a presymplectic structure on $M$; it degenerates
for example at the closed geodesics of $X$.

Despite this degeneracy the circle action is generated by a
Hamiltonian, namely the energy function $E\colon M\to\R$ given by
$E(\gamma)=\frac12\int_{S^1}\norm{d\gamma}^2$.  The
Duistermaat-Heckman theorem tells us to integrate the functional
$e^{-tE}$ times a ``Liouville'' volume form on $M$.  The
``Riemannian'' volume form on $M$ is the Wiener measure $d\gamma$, and
just as in the finite-dimensional case we must multiply this by the
(regularized) Pfaffian of the skew symmetric endomorphism of $TM$
defined by the presymplectic form.  This Pfaffian exists if the
manifold $X$ has a $\Spin$-structure, and the Duistermaat-Heckman
integral $\int_Me^{-tE(\gamma)}\Pf(\omega)\,d\gamma$ is seen to be the
index of the associated Dirac operator $\dirac$.  The
Duistermaat-Heckman theorem then says that this integral localizes to
the fixed point set $M^{S^1}$, which is a copy of $X$.  By calculating
the weights of the action on the normal bundle of $X$ in $M$ one
arrives at the A-roof genus and thus concludes that
$\index(\dirac)=\hat{A}(X)$.

\subsection*{Nitta's theorem}

Is there a generalization of the Duistermaat-Heckman theorem to
Poisson manifolds?  In general this seems too much to ask for, but a
reasonable compromise was found by Yasufumi Nitta
\cite{nitta;duistermaat-heckman}, \cite{nitta;reduction-calabi-yau}.
I state his result in a more general form obtained by (my student and
Hans' grand-student) Yi Lin~\cite{lin;equivariant-twisted}.  Let
$(M,\rho)$ be a compact \emph{generalized Calabi-Yau manifold}, that
is a $2n$-dimensional manifold equipped with a (possibly twisted)
generalized complex structure defined by a pure spinor
$\rho\in\Omega^*(M,\C)$ with the property that the $2n$-form
$\nu=(\rho,\bar{\rho})$ is a volume form.  (Such manifolds are
discussed in more detail in Gil Cavalcanti's lecture notes in these
proceedings.)  Let $T$ be a torus acting on $M$ in a Hamiltonian
fashion with generalized moment map $\Phi\colon M\to\t^*$.  Then the
pushforward measure $\Phi_*(\nu)$ on $\t^*$ is equal to a piecewise
polynomial function times Lebesgue measure.  For a symplectic manifold
$(M,\omega)$ we have $\rho=e^{i\omega}$ and $\nu=\omega^n/n!$, and so
this assertion is exactly the Fourier transformed version of the
classical Duistermaat-Heckman theorem.


\bibliographystyle{amsplain}

\def\cprime{$'$}
\providecommand{\bysame}{\leavevmode\hbox to3em{\hrulefill}\thinspace}
\providecommand{\MR}{\relax\ifhmode\unskip\space\fi MR }
\providecommand{\MRhref}[2]{%
  \href{http://www.ams.org/mathscinet-getitem?mr=#1}{#2}
}
\providecommand{\href}[2]{#2}


\end{document}